\documentclass{amsart}

\usepackage{fancyhdr}
\usepackage{lastpage}
\usepackage{stmaryrd,yhmath}

\newcommand{\bdis}{\begin{displaymath}}
\newcommand{\edis}{\end{displaymath}}
\newcommand{\be}{\begin{equation}}
\newcommand{\ee}{\end{equation}}

\newcommand{\mbb}{\mathbb}
\newcommand{\mcal}{\mathcal}

\newcommand{\vp}{\varphi}
\newcommand{\vt}{\vartheta}

\newcommand{\tit}{\tilde{t}}

\newcommand{\zf}{\zeta\left(\frac{1}{2}+it\right)}

\DeclareMathOperator*{\ssum}{\sum\sum}

\newtheorem{lemma}[]{Lemma}

\theoremstyle{definition}

\newtheorem{cor}[]{Corollary}

\theoremstyle{remark}
\newtheorem{remark}[]{Remark}

\newtheorem*{mydef11}{{\bf Theorem 1}}

\newtheorem*{mydef12}{{\bf Theorem 2}}

\newtheorem*{mydef13}{{\bf Theorem 3}}

\newtheorem*{mydef5A}{{\bf Lemma A}}
\newtheorem*{mydef5B}{{\bf Lemma B}}

\numberwithin{equation}{section}

\begin{document}

\title{On the order of the Titchmarsh's sum in the theory of the Riemann zeta-function and on the biquadratic effect in the information theory}

\author{Jan Moser}

\address{Department of Mathematical Analysis and Numerical Mathematics, Comenius University, Mlynska Dolina M105, 842 48 Bratislava, SLOVAKIA}

\email{jan.mozer@fmph.uniba.sk}


\begin{abstract}
We obtain in this paper the solution of the classical problem on the order of the Titchmarsh's sum (1934). Simultaneously, we obtain a connection of this problem and the Kotelnikoff-Whittaker-Nyquist's theorem
from the information theory.
\end{abstract}

\maketitle

\section{Introduction}

The following is the translation of the paper of reference \cite{9} into English. In this paper we obtain the solution of the classical problem on the order of the complicated Titchmarsh's sum

\bdis
\sum_{\nu=M+1}^N Z^2(t_\nu)Z^2(t_{\nu+1}) .
\edis

In connection with this we also obtain an analog of the biquadratic effect for

\bdis
Z(t),\ t\in [T,2T] .
\edis

It follows that the continuous signal defined by the function $Z(t)$ obeys the theorem of Kotelnikoff-Whittaker-Nyquist from the information theory. \\

Let us remind the definition of the Riemann zeta-function

\bdis
\zeta(s)=\prod_p \frac{1}{1-\frac{1}{p^s}},\ s=\sigma+it,\ \sigma>1
\edis
($p$ runs over through the set of all primes), and the analytic continuation of this function to all $s\in\mbb{C},\ s\not=1$. Riemann defined also the real-valued function

\be \label{1.1}
\begin{split}
& Z(t)=e^{i\vartheta(t)}\zf,\ \vt(t)=-\frac t2\ln\pi+\text{Im}\ln\Gamma\left(\frac 14+i\frac t2\right)= \\
& =\frac t2\ln\frac{t}{2\pi}-\frac{t}{2}-\frac{\pi}{8}+\mcal{O}\left(\frac 1t\right)
\end{split}
\ee
(see \cite{10}, (35), (44), (62), \cite{11}, p. 98). From this, it follows that the properties of the signal generated by the Riemann's function are connected with the law of the distribution of the primes in the series
of all positive integers, and this is to be regarded as a pleasant circumstance from the point of view of the Pythagorean philosophy of the Universe.

\section{The asymptotic formula for the Titchmarsh's sum}

\subsection{}

In 1934 Titchmarsh presented the following hypothesis (see \cite{11}, p. 105):\ there is $A>0$ such that
\bdis
\sum_{\nu=M+1}^N Z^2(t_\nu)Z^2(t_{\nu+1})=\mcal{O}(N\ln ^A N)
\edis
where $M$ is sufficiently big fixed number and $\{ t_\nu\}$ is the sequence defined by the condition (comp. \cite{11}, p. 99)
\be \label{2.1}
\vt(t_\nu)=\pi\nu ,\ \nu=1,2,\dots \ .
\ee
In 1980, we have proved this hypothesis with $A=4$ (see \cite{4}, (4)). The following estimate (see \cite{4}, (6))

\be \label{2.2}
\sum_{\tilde{t}_{M+1}\leq \tilde{t}_\nu\leq T} Z^4(\tilde{t}_\nu)=\mcal{O}(T\ln^5T)
\ee
was the key to our proof. The sequence $\{\tilde{t}_\nu\}$ is defined by the formula (see \cite{4}, (5))
\be \label{2.3}
\vt(\tilde{t}_\nu)=\frac{\pi}{2}\nu,\ \nu=1,2,\dots \ .
\ee

In 1983, we have improved the estimate (\ref{2.2}), namely, we have proved the asymptotic formula (see \cite{7})

\bdis
\sum_{\tilde{t}_{M+1}\leq \tilde{t}_\nu\leq T} Z^4(\tilde{t}_\nu)\sim\frac{1}{2\pi^3}T\ln^5T,\ T\to \infty .
\edis

In this paper we obtain the solution of the classical Titchmarsh's problem. In reality, we obtain the general autocorrelative formulae for the
function $Z^2(t)$, and from these formulae, as a special case, we obtain the desired result. Namely, the following main Theorem holds true.

\begin{mydef11}
\be \label{2.4}
\begin{split}
& \sum_{T\leq t_\nu\leq 2T}Z^2\{ t_\nu+k\rho_1(\nu)\}Z^2\{ t_\nu+l\rho_1(\nu)\}=\\
& =\left\{\begin{array}{lcl}\frac{3}{4\pi^5(k-l)^2}T\ln^5T+\mcal{O}(MT\ln^4T) & , & k\not=l \\
\frac{1}{4\pi^3}T\ln^5T+\mcal{O}\{(M+1)T\ln^4T\} & , & k=l \end{array} \right.,
\end{split}
\ee
where
\be \label{2.5}
\rho_1(\nu)=\frac{2\pi}{\ln\frac{t_\nu}{2\pi}},\ k,l=0,\pm 1,\pm 2,\dots ,\pm M,\ M=\mcal{O}(\psi) ,
\ee
and $\psi=\psi(T)$ is a function arbitrarily slowly increasing to $\infty$ as $T\to\infty$.
\end{mydef11}

\subsection{}

We obtain the final result on the order of Titchmarsh's sum from our Theorem 1 as follows. First of all we have (see (\ref{2.4}))

\be \label{2.6}
\sum_{T\leq t_\nu\leq 2T}Z^2(t_\nu)Z^2\{ t_\nu+\rho_1(\nu)\}=\frac{3}{4\pi^5}T\ln^5T+\mcal{O}(T\ln^4T) ,
\ee
$k=0,\ l=1;\ M=1$. Since (see \cite{3}, (42))

\be \label{2.7}
t_{\nu+1}-t_\nu=\frac{2\pi}{\ln\frac{t_\nu}{2\pi}}+\mcal{O}\left(\frac{1}{t_\nu\ln t_\nu}\right)=\rho_1(\nu)+
\mcal{O}\left(\frac{1}{t_\nu\ln^2 t_\nu}\right)
\ee

then we obtain by the usual estimates

\bdis
Z(t)=\mcal{O}(t^{1/6}\ln t),\ Z'(t)=\mcal{O}(t^{1/6}\ln ^2t)
\edis

the following

\be \label{2.8}
Z^2(t_{\nu+1})=Z^2\left\{t_\nu+\rho_1(\nu)+\mcal{O}\left(\frac{1}{t_\nu\ln^2 t_\nu}\right)\right\}=
Z^2\{ t_\nu+\rho_1(\nu)\}+\mcal{O}\left(\frac{\ln T}{T^{2/3}}\right) .
\ee

Next, we obtain by (\ref{2.8}) and by the formula

\bdis
\sum_{T\leq t_\nu\leq 2T} 1=\frac{1}{2\pi}T\ln T+\mcal{O}(T)
\edis

the following

\be \label{2.9}
\begin{split}
& \sum_{T\leq t_\nu\leq 2T}Z^2(t_\nu)Z^2(t_{\nu+1})=\sum_{T\leq t_\nu\leq 2T}Z^2(t_\nu)Z^2\{ t_\nu+\rho_1(\nu)\}+\\
& +\mcal{O}(T^{1/3}\ln^2T\cdot T^{-2/3}\ln T\cdot T\ln T)= \\
& =\sum_{T\leq t_\nu\leq 2T}Z^2(t_\nu)Z^2\{ t_\nu+\rho_1(\nu)\}+\mcal{O}(T^{2/3}\ln^4T) .
\end{split}
\ee

Hence, by (\ref{2.6}), (\ref{2.9}) we obtain

\begin{cor}
\be \label{2.10}
\sum_{T\leq t_\nu\leq 2T}Z^2(t_\nu)Z^2(t_{\nu+1})=\frac{3}{4\pi^5}T\ln^5T+\mcal{O}(T\ln^4T) .
\ee
\end{cor}

\begin{remark}
The order of the Titchmarsh's sum is determined by the asymptotic formula (\ref{2.10}).
\end{remark}

\section{Main lemmas and the conclusion of the proof of Theorem 1}

The following main lemmas hold true.

\begin{mydef5A}
\be \label{3.1}
\begin{split}
& \sum_{T\leq \tilde{t}_\nu\leq 2T}Z^2\{ \tilde{t}_\nu+k\rho_2(\nu)\}Z^2\{ \tilde{t}_\nu+l\rho_2(\nu)\}= \\
& =\left\{\begin{array}{lcl}\frac{3}{2\pi^5(k-l)^2}T\ln^5T+\mcal{O}(MT\ln^4T) & , & k\not=l \\
\frac{1}{2\pi^3}T\ln^5T+\mcal{O}\{(M+1)T\ln^4T\} & , & k=l \end{array} \right.,
\end{split}
\ee
where
\be \label{3.2}
\rho_2(\nu)=\frac{2\pi}{\ln\frac{\tilde{t}_\nu}{2\pi}} ,
\ee
and $k,l,M,\psi$ fulfills the conditions of Theorem 1.
\end{mydef5A}

\begin{mydef5B}
\be \label{3.3}
\sum_{T\leq \tilde{t}_\nu\leq 2T}(-1)^\nu Z^2\{ \tilde{t}_\nu+k\rho_2(\nu)\}Z^2\{\tilde{t}_\nu+l\rho_2(\nu)\}=
\mcal{O}\{(M+1)T\ln^4T\} .
\ee
\end{mydef5B}

Using these lemmas we easily conclude the proof of Theorem 1. Namely, by adding (\ref{3.1}) and (\ref{3.2}) we have

\be \label{3.4}
\begin{split}
& \sum_{T\leq \tilde{t}_\nu\leq 2T}Z^2\{\tilde{t}_{2\nu}+k\rho_2(2\nu)\}Z^2\{\tilde{t}_{2\nu}+l\rho_2(2\nu)\}= \\
& =\left\{\begin{array}{lcl}\frac{3}{4\pi^5(k-l)^2}T\ln^5T+\mcal{O}(MT\ln^4T) & , & k\not=l \\
\frac{1}{4\pi^3}T\ln^5T+\mcal{O}\{(M+1)T\ln^4T\} & , & k=l \end{array} \right. .
\end{split}
\ee

Since (see (\ref{2.1}), (\ref{2.3}), (\ref{2.5}), (\ref{3.2})), $\tilde{t}_{2\nu}=t_\nu,\ \rho_2(2\nu)=\rho_1(\nu)$, then
the formula (\ref{2.4}) follows from (\ref{3.4}). \\

The proofs of our lemmas A and B are situated in the parts 5-7 and 8-10, respectively.

\section{Biquadratic effect and the connection with the theorem of Kotelnikoff-Whittaker-Nyquist}

\subsection{}

Next, we obtain from (\ref{2.4}), $k=l=M=0$ the following

\begin{cor}
\be \label{4.1}
\sum_{T\leq t_\nu\leq 2T}Z^4(t_\nu)=\frac{1}{4\pi^3}T\ln^5T+\mcal{O}(T\ln^4T) .
\ee
\end{cor}

Since (see \cite{2}, p. 227, comp. \cite{12}, p. 125) we have

\be \label{4.2}
\int_T^{2T} Z^4(t){\rm d}t=\frac{1}{2\pi^2}T\ln^4T+\mcal{O}(T\ln^3T) ,
\ee
and by (\ref{4.1})
\be \label{4.3}
\frac{2\pi}{\ln T}\sum_{T\leq t_\nu\leq 2T}Z^4(t_\nu)=\frac{1}{2\pi^2}T\ln^4T+\mcal{O}(T\ln^3T) ,
\ee

then from (\ref{4.2}) by (\ref{4.3}) one obtains the following statement.

\begin{mydef12}
\be \label{4.4}
\int_T^{2T}Z^4(t){\rm d}t\sim \frac{2\pi}{\ln T}\sum_{T\leq t_\nu\leq 2T}Z^4(t_\nu),\ T\to\infty .
\ee
\end{mydef12}

We shall give another example of the relation of the type (\ref{4.4}). First of all, we have the Hardy-Littlewood mean-value theorem

\be \label{4.5}
\int_T^{T+U}Z^2(t){\rm d}t\sim U\ln T,\ U=\sqrt{T}\ln T,\ T\to\infty .
\ee

Next, we have proved a discrete analog of the formula (\ref{4.5}) (see \cite{5}, (6), comp. \cite{6}, (10); $H\to U,\ \tau^\prime=0$)

\be \label{4.6}
\sum_{T\leq t_\nu\leq 2T}Z^2(t_\nu)\sim \frac{1}{2\pi}U\ln^2T,\ T\to\infty .
\ee

The next statement follows immediately from the formulae (\ref{4.5}) and (\ref{4.6}).

\begin{mydef13}
\be \label{4.7}
\int_T^{T+U}Z^2(t){\rm d}t\sim\frac{2\pi}{\ln T}\sum_{T\leq t_\nu\leq T+U}Z^2(t_\nu),\ T\to\infty .
\ee
\end{mydef13}

\subsection{}

Let us remind some facts about the information theory. If we use the continuous signals in the information theory then the Kotelnikoff-Whittaker-Nyquist
theorem is the basic mathematical instrument. Namely, in the radioengineering the following \emph{empirical} rule is used (see \cite{1}, pp. 81, 86, 96,
97): if the length of the signal $F(t)$ is approximately $T$ (for example, $t\in [0,T]$), and the spectrum of the given signal $F(t)$ is bounded
approximately by the frequency $w$, and if $2Tw\gg 1$ then we have

\be \label{4.8}
F(t)\approx \sum_{n=0}^{2Tw}\frac{\sin(2\pi wt-n\pi)}{2\pi wt-n\pi}F\left(\frac{n}{2w}\right),\ t\in [0,T] ,
\ee

\be \label{4.9}
\int_0^TF^2(t){\rm d}t\approx \frac{1}{2w}\sum_{n=0}^{2wT}F^2\left(\frac{n}{2w}\right) .
\ee

The quantities
\bdis
\frac{1}{2w},\quad \int_0^TF^2(t){\rm d}t
\edis

are said to be the length of the Nyquist interval and the quadratic effect, respectively.

\begin{remark}
Since (see (\ref{2.7}))
\bdis
t_{\nu+1}-t_\nu\sim\frac{2\pi}{\ln T},\ t_\nu\in [T,T+U];\ [T,2T],\ T\to\infty
\edis
then by the asymptotic formulae (\ref{4.7}) and (\ref{4.4}) we have expressed the quadratic effect (comp. (\ref{4.9})) and the biquadratic effect,
respectively, of the signal defined by the function $Z(t)$. The following length of the Nyquist's interval
\bdis
\frac{1}{2w}\sim \frac{2\pi}{\ln T} ,\ T\to \infty
\edis
corresponds with these effects.
\end{remark}

\begin{remark}
We have proved also the analog of the equation (\ref{4.8}) for the function $Z(t)$, in the sense of the discrete mean-square.
\end{remark}

\section{The formula for $Z^2\{\tilde{t}_\nu+k\rho_2(\nu)\}Z^2\{\tilde{t}_\nu+l\rho_2(\nu)\}$}

We use the Hardy-Littlewood formula (comp. \cite{4}, (24))

\bdis
Z^2(t)=2\sum_{n\leq t_1}\frac{d(n)}{\sqrt{n}}\cos\{ 2\vt(t)-t\ln n\}+\mcal{O}(\ln T),\ t_1=\frac{t}{2\pi} , \ t\in [T,2T]
\edis

where $d(n)$ is the number of divisors of $n$. Then we have

\be \label{5.1}
\begin{split}
& Z^2\{\tilde{t}_\nu+k\rho_2(\nu)\}=2\sum_{n\leq t_2}\frac{d(n)}{\sqrt{n}}\cos\{ 2\vt(\tilde{t}+k\rho_2(\nu))-\tilde{t}_\nu\ln n-k\rho_2(\nu)\ln n\}+ \\
& + \mcal{O}(\ln T),\ t_2=\frac{\tilde{t}_\nu}{2\pi},\ \tilde{t}_\nu\in [T,2T] ,
\end{split}
\ee

where the inequality $n\leq t_2$ follows from the estimate

\bdis
\frac{\tit_\nu+k\rho_2(\nu)}{2\pi}-t_2=\frac{k\rho_2(\nu)}{2\pi}=\mcal{O}\left(\frac{M+1}{\ln T}\right)=o(1) ,
\edis

(see (\ref{2.5}), (\ref{3.2})). Since (see \cite{12},. p. 221)

\be \label{5.2}
\vt'(t)=\frac 12\ln\frac{t}{2\pi}+\mcal{O}\left(\frac{1}{t}\right),\quad \vt''(t)\sim \frac{1}{2t} ,
\ee

then (see (\ref{2.3}))

\be \label{5.3}
2\vt(\tit_\nu+k\rho_2(\nu))=\pi\nu+2k\pi+\mcal{O}\left(\frac{M+1}{T\ln T}\right),\ \tit_\nu\in [T,2T] .
\ee

Since the remainder in (\ref{5.3}) generates the error

\bdis
\mcal{O}\left(\frac{1}{T^{1/2-\epsilon}\ln T}\right)
\edis

($\epsilon>0$ is arbitrarily small) in (\ref{5.1}), then we obtain by (\ref{5.3}) the following

\be \label{5.4}
\begin{split}
& Z^2\{\tilde{t}_\nu+k\rho_2(\nu)\} = \\
& =2(-1)^\nu\sum_{n\leq t_2}\frac{d(n)}{\sqrt{n}}\cos\{\tit_\nu\ln n+k\rho_2(\nu)\ln n\}+\mcal{O}(\ln T),\ \tit_\nu\in [T,2T] .
\end{split}
\ee

Consequently

\be \label{5.5}
Z^2\{\tilde{t}_\nu+k\rho_2(\nu)\}Z^2\{\tilde{t}_\nu+l\rho_2(\nu)\}=S+R_1+R_2+\mcal{O}(\ln^2 T)
\ee

where

\be \label{5.6}
\begin{split}
& S=2\sum_{n\leq t_2}\frac{d^2(n)}{n}\cos\{(k-l)\rho_2(\nu)\ln n\}+ \\
& + 2\ssum_{m,n\leq t_2,m\not=n} \frac{d(m)d(n)}{\sqrt{mn}}\cos\{\tit_\nu\ln\frac mn + k\rho_2(\nu)\ln m-l\rho_2(\nu)\ln n\}+ \\
& + 2\sum_{n\leq t_2}\frac{d^2(n)}{n}\cos\{ 2\tit_\nu\ln n+(k+l)\rho_2(\nu)\ln n\}+ \\
& + 2\ssum_{m,n\leq t_2,m\not=n}\frac{d(m)d(n)}{\sqrt{mn}}\cos\{\tit_\nu\ln (mn)+k\rho_2(\nu)\ln m+l\rho_2(\nu)\ln n\}= \\
& = S_1+S_2+S_3+S_4 ,
\end{split}
\ee

and (see (\ref{5.4}))

\be \label{5.7}
\begin{split}
& R_1=\mcal{O}\left(\ln T\left|\sum_{n\leq t_2}\frac{d(n)}{\sqrt{n}}\cos\{\tit_\nu\ln n+k\rho_2(\nu)\ln n\}\right|\right)= \\
& = \mcal{O}\left(Z^2\{\tit_\nu+k\rho_2(\nu)\}\ln T\right)+\mcal{O}(\ln^2T) , \\
& R_2=\mcal{O}\left(Z^2\{\tit_\nu+l\rho_2(\nu)\}\ln T\right)+\mcal{O}(\ln^2T) .
\end{split}
\ee

\section{The main term in the asymptotic formula (\ref{3.1})}

The following lemma holds true.

\begin{lemma}
\be \label{6.1}
\sum_{T\leq \tit_\nu\leq 2T} S_1
=\left\{\begin{array}{lcl}\frac{3}{2\pi^5(k-l)^2}T\ln^5T+\mcal{O}(T\ln^4T) & , & k\not=l \\
\frac{1}{2\pi^3}T\ln^5T+\mcal{O}\{T\ln^4T\} & , & k=l \end{array} \right. .
\ee
\end{lemma}

\begin{proof}
We have (see (\ref{5.6}))
\be \label{6.2}
\begin{split}
& S_1=2S_{11} , \\
& S_{11}=\sum_{n\leq t_2}\frac{d^2(n)}{n}\cos(\alpha \ln n),\ t_2=\frac{\tit_\nu}{2\pi},\ \alpha=(k-l)\rho_2(\nu) .
\end{split}
\ee
\begin{itemize}
\item[(A)] Let $k\not=l$, i.e. $\alpha\not=0$. \\

By using of the partial summation and the Ramanujan's formula (see \cite{2}, p. 296)
\bdis
D(x)=\sum_{n=1}^x d^2(n)=\frac{1}{\pi^2}x\ln^3 x+\mcal{O}(\ln^2x),\ x=\frac{\tit_\nu}{2\pi} ,
\edis
we obtain ($D(0)=0$)
\be \label{6.3}
\begin{split}
& S_{11}=\sum_{n=1}^{[x]}\{ D(n)-D(n-1)\}\frac 1n\cos(\alpha\ln n)= \\
& = \sum_{n=1}^{[x]}D(n)\left\{\frac{\cos(\alpha\ln n)}{n}-\frac{\cos(\alpha\ln(n+1))}{n+1}\right\}+\mcal{O}(\ln^3x)=\\
& = \sum_{n=1}^{[x]}D(n)\int_n^{n+1}\{\cos(\alpha\ln v)+\alpha\sin(\alpha\ln v)\}\frac{{\rm d}v}{v^2}+\mcal{O}(\ln^3x)= \\
& = \frac{1}{\pi^2}\int_1^x\ln^3v\{\cos(\alpha\ln v)+\alpha\sin(\alpha\ln v)\}\frac{{\rm d}v}{v}+\mcal{O}(\ln^3x) = \\
& = \frac{1}{\pi^2}\int_0^{\ln x}\{ w^3\cos(\alpha w)+\alpha w^3\sin(\alpha w)\}{\rm d}w+\mcal{O}(\ln^3 x) = \\
& = \frac{1}{\pi^2}F(x,\alpha)+\mcal{O}(\ln^3x) .
\end{split}
\ee
Next, by using of the simple integration by parts, we obtain
\be \label{6.4}
\begin{split}
& F(x,\alpha)=\left[\left(\frac{3w^2}{\alpha^2}-\frac{6}{\alpha^4}+\frac{6w}{\alpha^2}-w^3\right)\cos(\alpha w)+ \right. \\
& \left. +\left(\frac{w}{\alpha}-\frac{6w}{\alpha^3}+\frac{3w^2}{\alpha}-\frac{6}{\alpha^3}\right)\sin(\alpha w) \right]_{0}^{\ln x} = \\
& =\left(3\frac{\ln^2 x}{\alpha^2}-\frac{6}{\alpha^4}\right)\cos(\alpha \ln x)+\left(\frac{\ln^3x}{\alpha}-6\frac{\ln x}{\alpha^3}\right)
\sin(\alpha\ln x)+ \\
& + \frac{6}{\alpha^4}+\mcal{O}(\ln^3T)=\frac{3}{4\pi^2(k-l)^2}\ln^4\frac{\tit_\nu}{2\pi}+\mcal{O}(\ln^3T) ,
\end{split}
\ee
since (see (\ref{3.2}), (\ref{6.2}))
\bdis
\alpha\ln x=2\pi (k-l) .
\edis
Consequently (see (\ref{6.2})-(\ref{6.4}))
\be \label{6.5}
S_1=\frac{3}{2\pi^4(k-l)^2}\ln^4T+\mcal{O}(\ln^3T) ,
\ee
and, of course,
\be \label{6.6}
\sum_{T\leq\tit_\nu\leq 2T}1=\frac 1\pi T\ln T+\mcal{O}(T) .
\ee
Hence, by (\ref{6.5}), (\ref{6.6}) the first formula in (\ref{6.1}) follows.

\item[(B)] Let $k=l$, i.e. $\alpha=0$. \\

Putting $\alpha=0$ in the fifth line of the formula (\ref{6.3}), we obtain
\bdis
S_{11}=\frac{1}{4\pi^2}\ln^4\frac{T}{2\pi}+\mcal{O}(\ln^3T) .
\edis
This, together with (\ref{6.6}) gives the second formula in (\ref{6.1}) .
\end{itemize}
\end{proof}

\section{The estimates of the remaining terms}

\subsection{}

The following lemma holds true.

\begin{lemma}
\be \label{7.1}
\sum_{T\leq\tit_\nu\leq 2T}S_2=\mcal{O}(T\ln^4T) .
\ee
\end{lemma}

\begin{proof}
Let
\be \label{7.2}
\tit_k\leq 2T\leq \tit_{k+1},\ \tau=\max\{T,2\pi n,2\pi m\} .
\ee
We have (see (\ref{5.6})

\be \label{7.3}
\begin{split}
& \sum_{T\leq\tit_\nu\leq 2T}S_2=2\ssum_{m,n\leq \tit_k/2\pi, m\not=n} \frac{d(m)d(n)}{\sqrt{mn}}U_2 , \\
& U_2=\sum_{\tau\leq\tit_\nu\leq 2T}\cos\{\tit_\nu\ln\frac mn +h_1(\nu)\},\ h_1(\nu)=k\rho_2(\nu)\ln m-l\rho_2(\nu)\ln n .
\end{split}
\ee
First of all we put
\be \label{7.4}
\begin{split}
& U_2=U_{21}-U_{22} , \\
& U_{21}=\sum_{\tau\leq \tit_\nu\leq 2T}\cos\{h_1(\nu)\}\cos\left\{\tit_\nu\ln\frac mn\right\} , \\
& U_{21}=\sum_{\tau\leq \tit_\nu\leq 2T}\sin\{h_1(\nu)\}\sin\left\{\tit_\nu\ln\frac mn\right\} .
\end{split}
\ee

It is sufficient to estimate the term $U_{21}$. Since
\bdis
h_1(\nu)=\mcal{O}(M),\ \tit_\nu\in [T,2T],
\edis
then
\bdis
h_1(\nu)\in [-AM,AM] .
\edis
Now, we divide the segment $[-AM,AM]$ on $\mcal{O}(M)$ parts in such a way that on each part of our segment the following is true: either
\bdis
0\leq \cos\{h_1(\nu)\}\leq 1 ,
\edis
or
\bdis
0\leq -\cos\{h_1(\nu)\}\leq 1 ,
\edis
and the sequences
\bdis
\cos\{h_1(\nu)\},-\cos\{h_1(\nu)\}
\edis
are monotone. If we use the Abel's transformation on every of those parts, we obtain
\be \label{7.5}
|U_{21}|\leq AM\cdot\max_{\tau_1,\tau_2;\ \tau\leq \tau_1<\tau_2\leq 2T}
\left|\sum_{\tau_1\leq\tit_\nu\leq \tau_2}\cos\left\{\tit_\nu\ln\frac mn\right\}\right| ,
\ee
Of course, instead of the sum in (\ref{7.5}) we may estimate the following sum
\bdis
U_{211}=\sum_{\tau\leq \tit_\nu\leq \tau_1\leq 2T} \cos\left\{\tit_\nu\ln\frac mn\right\} ,
\edis
and for this sum the method explained in \cite{4}, pp. (30)-(37) is applicable. We then obtain from (\ref{7.5}) the estimate
\bdis
U_{21}=\mcal{O}\left(\frac{(M+1)\ln T}{\left|\ln\frac mn\right|}\right) .
\edis
This estimate is valid also for $U_{22}$ and, by (\ref{7.4}), for $U_2$. Hence, from (\ref{7.3}) (comp. \cite{4}, (38), (39)) we obtain (\ref{7.1}).
\end{proof}

\subsection{}

On the basis of \cite{4}, (40)-(59), \cite{7}, (5), (6), via a similar way, we obtain

\begin{lemma}
\be \label{7.6}
\sum_{T\leq\tit_\nu\leq 2T}S_4=\mcal{O}(T\ln^4T) .
\ee
\end{lemma}

Next, on the basis \cite{7}, (8)-(12), we obtain

\begin{lemma}
\be \label{7.7}
\sum_{T\leq\tit_\nu\leq 2T}S_3=\mcal{O}(T\ln T) .
\ee
\end{lemma}

Consequently, by (\ref{5.6}), (\ref{7.1}), (\ref{7.6}), (\ref{7.7}) we have

\be \label{7.8}
\sum_{T\leq\tit_\nu\leq 2T}S=\sum_{T\leq\tit_\nu\leq 2T}S_1+ \mcal{O}(T\ln^4T) .
\ee

\subsection{}

From the Riemann-Siegel formula

\bdis
Z(t)=2\sum_{n\leq t_3}\frac{1}{\sqrt{n}}\cos\{\vt(t)-t\ln n)+\mcal{O}(t^{-1/4}),\ t_3=\sqrt{\frac{t}{2\pi}}
\edis

we easily obtain the estimate

\be \label{7.9}
\sum_{T\leq\tit_\nu\leq 2T}Z^2\{\tit_\nu+k\rho_2(\nu)\}=\mcal{O}(T\ln^2T) .
\ee

Thus (see (\ref{5.7}), (\ref{6.6}), (\ref{7.9})) we have

\be \label{7.10}
\sum_{T\leq\tit_\nu\leq 2T}\{ R_1+R_2+\mcal{O}(\ln^2T)\}=\mcal{O}(T\ln^3T) .
\ee

Finally, from (\ref{5.5}) by (\ref{6.1}), (\ref{7.8}), (\ref{7.10}) we obtain (\ref{3.1}).

\section{The formula for $(-1)^\nu Z^2\{\tit_\nu+k\rho_2(\nu)\}Z^2\{\tit_\nu+l\rho_2(\nu)\}$}

First of all (see (\ref{5.5}))

\be \label{8.1}
\begin{split}
& (-1)^\nu Z^2\{\tit_\nu+k\rho_2(\nu)\}Z^2\{\tit_\nu+l\rho_2(\nu)\}= \\
& = \bar{S}+(-1)^\nu (R_1+R_2)+\mcal{O}(\ln^2T) ,
\end{split}
\ee

where (see (\ref{5.6}))

\be \label{8.2}
\bar{S}=(-1)^\nu S_1+(-1)^\nu S_2+(-1)^\nu S_3+(-1)^\nu S_4 .
\ee

From (\ref{6.5}) we obtain

\be \label{8.3}
\sum_{T\leq\tit_\nu\leq 2T}(-1)^\nu S_1=\mcal{O}(T\ln^4T) .
\ee

Next we have (see (\ref{5.6}), (\ref{7.3}))

\bdis
\begin{split}
& \sum_{T\leq\tit_\nu\leq 2T}(-1)^\nu S_3=2\sum_{n\leq \tit_\nu/2\pi}\frac{d^2(n)}{n}\bar{U}_3 , \\
& \bar{U}_3=\sum_{\tau\leq\tit_\nu\leq 2T}\cos\{\pi\nu+2\tit_\nu\ln n+h_2(\nu)\} , \\
& h_2(\nu)=(k+l)\rho_2(\nu)\ln n .
\end{split}
\edis

The estimate of the sum $\bar{U}_3$ we may carry forward to the estimates of the following sums

\bdis
\bar{U}_{311}(r)=\sum_{\tau\leq\tit_\nu\leq \tau_1\leq 2T}\cos\left\{\pi\nu+2\tit_\nu\ln n-\frac{\pi}{2}r\right\},\ r=0,1
\edis

(comp. the part 7.1). We obtain the estimates of these sums by the van der Corput's lemma with the second derivative (see \cite{12}, p. 61). Consequently (comp. \cite{7}, (7)-(12)), we obtain
the estimate

\be \label{8.4}
\sum_{T\leq\tit_\nu\leq 2T}(-1)^\nu S_3=\mcal{O}(T\ln T) .
\ee

Next, we have (see (\ref{5.7}), (\ref{7.10}), (\ref{8.1}))

\be \label{8.5}
\sum_{T\leq\tit_\nu\leq 2T}\{ (-1)^\nu(R_1+R_2)+\mcal{O}(\ln^2T)\}=\mcal{O}(T\ln^3T) .
\ee

Now, the proof of the Lemma B lies on the two following lemmas.

\begin{lemma}
\be \label{8.6}
\sum_{T\leq\tit_\nu\leq 2T} (-1)^\nu S_2=\mcal{O}(T\ln^{7/2}T) .
\ee
\end{lemma}

\begin{lemma}
\be \label{8.7}
\sum_{T\leq\tit_\nu\leq 2T} (-1)^\nu S_4=\mcal{O}(T\ln^{3}T) .
\ee
\end{lemma}

\section{Proof of the Lemma 5}

We have (see (\ref{5.6}), (\ref{7.2}), (\ref{7.3}))

\be \label{9.1}
\begin{split}
& W_2=\sum_{T\leq\tit_\nu\leq 2T}(-1)^\nu S_2=2\ssum_{m,n\leq \tit_\nu/2\pi, m\not=n}\frac{d(m)d(n)}{\sqrt{mn}}\bar{U}_2 , \\
& \bar{U}_2=\sum_{T\leq\tit_\nu\leq 2T}\cos\left\{\pi\nu-t_\nu\ln\frac mn -h_1(\nu)\right\} .
\end{split}
\ee

Let us remind that the sequence $\{ g_\nu\}$ is defined by the formula (see \cite{8}, (6))

\be \label{9.2}
\vt_1(g_\nu)=\frac{\pi}{2}\nu,\ \nu=1,2,\dots
\ee

where

\be \label{9.3}
\begin{split}
& \vt_1(t)=\frac t2\ln\frac{t}{2\pi}-\frac t2-\frac{\pi}{8},\ \vt'_1(t)=\frac 12\ln\frac{t}{2\pi},\ \vt_1''(t)=\frac{1}{2t} ,
\end{split}
\ee

and (see (\ref{1.1}))

\be \label{9.4}
\vt(t)=\vt_1(t)+\mcal{O}\left(\frac 1t\right) .
\ee

It is clear (comp. (\ref{5.2}) with (\ref{9.3})) that the sequence $\{ g_\nu\}$ is more advisable for the estimation of the sum $U_2$. \\

Since (see (\ref{2.3}), (\ref{9.2}), (\ref{9.3}))

\bdis
\mcal{O}\left(\frac{1}{\tit_\nu}\right)=\vt_1(\tit_\nu)-\vt_1(g_\nu)=(\tit_\nu-g_\nu)\vt_1'(T_\nu),\ T_\nu\in (g_\nu,\tit_\nu); (\tit_\mu,g_\nu)
\edis

then (see (\ref{9.3}))

\be \label{9.5}
\tit_\nu-g_\nu=\mcal{O}\left(\frac{1}{T\ln T}\right),\ \tit_\nu\in [T,2T] .
\ee

Thus we obtain from (\ref{9.1}) by (\ref{6.6}), (\ref{9.5})

\be \label{9.6}
\begin{split}
& \bar{U}_2=\bar{U}_{21}+\mcal{O}(\ln T),\\
& \bar{U}_{21}=\sum_{\tau\leq g_\nu\leq 2T}\cos\left\{\pi\nu-g_\nu\ln\frac mn-h_1(\nu)\right\} .
\end{split}
\ee

The estimation of the sum $\bar{U}_{21}$ may be carry forward to the estimation of the following sums (comp. the part 7.1).

\be \label{9.7}
\bar{U}_{211}(r)=\sum_{\tau\leq g_\nu\leq\tau_1\leq 2T}\cos\{ 2\pi\Phi_2(\nu)\},\ r=0,1
\ee

where

\be \label{9.8}
\Phi_2(\nu)=\frac{\nu}{2}-\frac{g_\nu}{2\pi}\ln\frac mn - \frac r2 .
\ee

Let $m>n$. Since (see (\ref{9.2}))

\be \label{9.9}
\frac{{\rm d}g_\nu}{{\rm d}\nu}=\frac{\pi}{2\vt_1'(g_\nu)}=\frac{\pi}{\ln\frac{g_\nu}{2\pi}}
\ee

(here and in similar cases we assume that $g_\nu$ is defined by (\ref{9.2}) for all $\nu\geq 1$) then
\be \label{9.10}
\Phi_2'(\nu)=\frac 12-\frac{1}{2\ln\frac{g_\nu}{2\pi}}\ln\frac mn,\ \Phi_2''(\nu)>0 .
\ee

Because

\bdis
0<\frac{1}{2\ln\frac{g_\nu}{2\pi}}\ln\frac mn\leq \frac{\ln m}{2\ln\frac{T}{2\pi}}\leq \frac{\ln\frac T\pi}{2\ln\frac{T}{2\pi}}<\frac 12+\epsilon \ \Rightarrow \ |\Phi_2'(\nu)|<\frac 12 ,
\edis

then (see \cite{12}, p. 65, Lemma 4.8)

\bdis
\bar{U}_{211}(r)=\int_{\tau\leq g_\nu\leq\tau_1}\cos\{ 2\pi\Phi_2(\nu)\}{\rm d}\nu+\mcal{O}(1) .
\edis

Next, since $\Phi_2'(\nu)$ is increasing (see (\ref{9.10})) then

\bdis
\begin{split}
& \Phi_2'(\nu)\geq \frac 12-\frac{1}{2\ln\frac{\tau}{2\pi}}\ln\frac mn = \frac{1}{2\ln\frac{\tau}{2\pi}}\ln\frac{\tau n}{2\pi m};\ g_\nu\in [\tau,\tau_1] .
\end{split}
\edis

Let $n\geq 2$. If $2\pi m>T$ then $\tau=2\pi m$ (see (\ref{7.2})), and

\bdis
\Phi_2'(\nu)\geq \frac{\ln 2}{2\ln\tau}>\frac{A}{\ln T} .
\edis

If $T\geq 2\pi n$ then $\tau=T$, and

\bdis
\Phi_2'(\nu)\geq \frac{1}{2\ln\frac{T}{2\pi}}\ln\frac{Tn}{2\pi m}\geq \frac{\ln 2}{2\ln\frac{T}{2\pi}}>\frac{A}{\ln T} .
\edis

Thus, in the case $m>n\geq 2$ one obtains the following estimate (by \cite{12}, p. 61, Lemma 4.2)

\bdis
\bar{U}_{211}(r)=\mcal{O}(\ln T) .
\edis

Consequently, for $\bar{U}_2$ (see (\ref{9.6}), comp. (\ref{7.5})), we have

\bdis
\bar{U}_2=\mcal{O}((M+1)\ln T) ,
\edis

and for the corresponding part of $W_{21}$ ($m>n\geq 2$) of the sum $W_2$ (see (\ref{9.1})) we obtain

\be \label{9.11}
W_{21}=\mcal{O}\left\{(M+1)\ln T\ssum_{m,n\leq T} \frac{d(m)d(n)}{\sqrt{mn}}\right\}=\mcal{O}\{(M+1)T\ln^3T\}
\ee

since (see \cite{2}, pp. 297, 298)

\be \label{9.12}
\ssum_{m,n\leq x}\frac{d(m)d(n)}{\sqrt{mn}}=\mcal{O}\{ x\ln^2x\} .
\ee

Let $n=1$. Since (see (\ref{9.9}), (\ref{9.10}); $m\geq 2$)

\bdis
\Phi_2'(\nu)=\frac{\pi}{2}\frac{\ln m}{g_\nu\ln^3\frac{g_\nu}{2\pi}} > \frac{A}{T\ln^3T} ,
\edis

then we obtain, by the lemma with the second derivative,

\bdis
\bar{U}_{211}(r;n=1)=\mcal{O}(\sqrt{T}\ln^{3/2}T) ,
\edis

i.e. (see (9.6))
\bdis
\bar{U}_2(n=1)=\mcal{O}\{(M+1)\sqrt{T}\ln^{3/2}T\} .
\edis

Consequently, for the corresponding part $W_{22}$ of the sum $W_2$ (see (9.1)) we obtain

\be \label{9.13}
\begin{split}
& W_{22}=\mcal{O}\left\{(M+1)\sqrt{T}\ln^{3/2}T\sum_{n\leq T}\frac{d(n)}{\sqrt{n}}\right\}= \\
& = \mcal{O}\left\{(M+1)\sqrt{T}\ln^{3/2}T\sqrt{T}\left(\sum_{n\leq T}\frac{d^2(n)}{n}\right)^{1/2}\right\} = \\
& = \mcal{O}\{(M+1)T\ln^{7/2}T\}
\end{split}
\ee

because (see \cite{2}, p. 296)

\bdis
\sum_{n\leq x}\frac{d^2(n)}{n}=\frac{1}{4\pi^2}\ln^4x+\mcal{O}(\ln^3x) .
\edis

Thus in the case $m>n$ we have (see (9.11), (9.13))

\be \label{9.14}
W_2(m>n)=\mcal{O}\{(M+1)T\ln^{7/2}T\} .
\ee

Now, let $n>m$. In this case we have (see (9.8))

\bdis
\begin{split}
& 2\pi\Phi_2(\nu)=\pi\nu-g_\nu\ln\frac mn-\frac{\pi}{2}r=\pi\nu+g_\nu\ln\frac nm-\frac{\pi}{2}r= \\
& = 2\pi\nu-2\pi\left(\frac{\nu}{2}-\frac{g_\nu}{2\pi}\ln\frac nm+\frac r4\right)=2\pi\nu-2\pi\tilde{\Phi}_2 ,
\end{split}
\edis

i.e. in this case the following estimate

\be \label{9.15}
W_2(n>m)=\mcal{O}\{(M+1)T\ln^{7/2}T\}
\ee

holds true. Finally, from (9.1) by (9.14), (9.15) the estimate (8.6) follows.

\section{Proof of the Lemma 6}

We have (see (5.6), (7.2))

\be \label{10.1}
\begin{split}
& W_4=\sum_{T\leq\tit_\nu\leq 2T}(-1)^\nu S_4=2\ssum_{m,n\leq\tit_k/2\pi, m\not=n}\frac{d(m)d(n)}{\sqrt{mn}}\bar{U}_4 , \\
& \bar{U}_4=\sum_{\tau\leq\tit_\nu\leq 2T}\cos\{ \pi\nu-\tit_\nu\ln(mn)-h_3(\nu)\} , \\
& h_3(\nu)=k\rho_2(\nu)\ln m+l\rho_2(\nu)\ln n .
\end{split}
\ee

The estimate for the sum $\bar{U}_4$ may be carry forward to estimate the following sums (similarly to the case (9.6), (9.7))

\bdis
\bar{U}_{411}(r)=\sum_{\tau\leq g_\nu\leq \tau_1\leq 2T} \cos\{ 2\pi\Phi_4(\nu)\}+\mcal{O}(\ln T),\ r=0,1
\edis

where

\bdis
\Phi_4(\nu)=\frac{\nu}{2}-\frac{g_\nu}{2\pi}\ln(nm)+\frac r4 .
\edis

First of all, we have

\be \label{10.2}
\Phi_4'(\nu)=\frac 12-\frac{\ln(nm)}{2\ln\frac{g_\nu}{2\pi}} ,\quad \Phi_4''(\nu)>0 .
\ee

Since

\bdis
0<\frac{\ln(mn)}{2\ln\frac{g_\nu}{2\pi}}\leq \frac{\ln\left(\frac T\pi\right)^2}{2\ln\frac{T}{2\pi}}<1+\epsilon \ \Rightarrow \
|\Phi_4'(\nu)|\leq \frac 12+\epsilon
\edis

then

\be \label{10.3}
\bar{U}_{411}(r)=\int_{\tau\leq g_\nu\leq \tau_1}\cos\{ 2\pi\Phi_4(\nu)\}{\rm d}\nu+\mcal{O}(\ln T) .
\ee

\subsection{}

Let

\bdis
mn<\frac{T}{2\pi} .
\edis

Then $\tau=T$ and, since $\Phi_4'$ is increasing (see (10.2)), we have

\bdis
\Phi_4'(\nu)\geq \frac 12-\frac{\ln(mn)}{2\ln\frac{T}{2\pi}}=\frac{1}{2\ln Q_1}\ln\frac{Q_1}{mn}>0,\ Q_1=\frac{T}{2\pi} ,
\edis

\bdis
m,n<Q_1,\ g_\nu\in [T,\tau_1] .
\edis

Thus, by the lemma with the first derivative, we obtain (see (10.3))

\be \label{10.4}
\bar{U}_{411}(r;mn<Q_1)=\mcal{O}\left(\frac{\ln T}{\ln\frac{Q_1}{mn}}\right) .
\ee

Let

\bdis
\frac{T}{2\pi}-\alpha_1<\frac{T}{2\pi} ,
\edis

where $\alpha_1>0$ is a convenient number. Then we have the following contribution $\bar{W}$ into the sum $W_4$

\bdis
\bar{W}=\mcal{O}\left(\frac{T^{2\epsilon}}{\sqrt{T}}\alpha_1 T^{\epsilon}T\ln T\right)=\mcal{O}(T^{1/2+4\epsilon}),\
\frac{T}{2\pi}-\alpha_1\leq mn<\frac{T}{2\pi} .
\edis

\begin{remark}
If we assume that $Q_1\in\mbb{N}$ (we use this assumption also in other similar cases) then

\be \label{10.5}
W_{41}=W_4(mn<T/2\pi)=\mcal{O}(T^{1/2+4\epsilon})+\bar{W}_{41}(mn<Q_1) .
\ee
\end{remark}

We have (see (10.1), (10.4))

\be\label{10.6}
\begin{split}
& \bar{W}_{41}=\mcal{O}\left\{(M+1)\ln T\ssum_{mn<Q_1}\frac{d(m)d(n)}{\sqrt{mn}\ln\frac{Q_1}{mn}}\right\}=  \\
& = \mcal{O}\left\{(M+1)T^{2\epsilon}\ln T\ssum_{mn<Q_1}\frac{1}{\sqrt{mn}\ln\frac{Q_1}{mn}}\right\} = \\
& = \mcal{O}\left\{(M+1)T^{3\epsilon}\ln T\sum_{q<Q_1}\frac{1}{\sqrt{q}\ln\frac{Q_1}{q}}\right\} = \\
& = \mcal{O}\{(M+1)T^{3\epsilon}\ln T\sqrt{Q_1}\ln Q_1\}= \\
& = \mcal{O}\{(M+1)T^{1/2+4\epsilon}\} .
\end{split}
\ee

Hence (see (10.5), (10.6))

\be \label{10.7}
W_{41}=\mcal{O}\{(M+1)T^{1/2+4\epsilon}\} .
\ee

\subsection{}

Let

\bdis
mn>\frac{T}{\pi} .
\edis

Since $g_\nu\leq 2T$ then (see (10.2))

\bdis
-\Phi_4'(\nu)=\frac{\ln(mn)}{2\ln\frac{g_\nu}{2\pi}}-\frac 12>0 ,
\edis

i.e. $\{-\Phi_4'(\nu)\}$ is decreasing. Consequently,

\bdis
-\Phi_4'(\nu)\geq \frac{\ln(mn)}{\ln\frac{T}{\pi}}-\frac 12=\frac{1}{2\ln Q_2}\ln\frac{mn}{Q_2},\ Q_2=\frac{T}{\pi} ,
\edis

and (comp. (10.4))

\bdis
\bar{U}_{411}(r;mn>T/2\pi)=\mcal{O}\left(\frac{\ln T}{\ln\frac{mn}{Q_2}}\right),\ mn>Q_2 .
\edis

Let

\be \label{10.8}
W_{42}=W_4(mn>T/\pi)=W_{421}(Q_2<mn<2Q_2)+W_{422}(2Q_2\leq mn) .
\ee

First of all we have (similarly to the part 10.1)

\be \label{10.9}
W_{421}=\mcal{O}\{(M+1)T^{1/2+4\epsilon}\}
\ee

here the following known estimate was used

\bdis
\sum_{Q_2<q<2Q_2}\frac{1}{\sqrt{q}\ln\frac{q}{Q_2}}=\mcal{O}(\sqrt{Q_2}\ln Q_2) .
\edis

Next, since

\bdis
\ln\frac{mn}{Q_2}\geq \ln 2;\ mn\geq 2Q_2 ,
\edis

then we obtain (see (9.12))

\be \label{10.10}
\begin{split}
& W_{422}=\mcal{O}\left\{(M+1)\ln T\ssum_{m,n<T/\pi}\frac{d(m)d(n)}{\sqrt{mn}}\right\}= \\
& = \mcal{O}\{(M+1)T\ln^3T\}.
\end{split}
\ee

Hence, (see (10.8)-(10.10)) we have

\be \label{10.11}
W_{42}=\mcal{O}\{(M+1)T\ln^3T\} .
\ee

\subsection{}

Let

\bdis
\frac{T}{2\pi}\leq mn\leq \frac{T}{\pi} .
\edis

Since $g_\nu\in [T,2T]$ then, in the case (10.12), the function $\Phi_4'(\nu)$ (see (10.2)) has only one zero $\bar{\nu}$ (since $\Phi_4'$ is
increasing).

\subsubsection{}

If $\tau\leq g_\nu\leq \tau_1$ then (see (10.3))

\be \label{10.13}
\begin{split}
& \bar{U}_{411}(r)=\int_{\tau\leq g_\nu\leq g_{\bar{\nu}}-A_1}+\int_{g_{\bar{\nu}}+A_2\leq g_\nu\leq \tau_1}+\mcal{O}(1)+\mcal{O}(\ln T) = \\
& = \bar{U}^1_{411}(r)+\bar{U}^2_{411}(r)+\mcal{O}(\ln T)
\end{split}
\ee

with evident adaptations if $g_{\bar{\nu}}=\tau,\tau_1$ ($0<A_1<A_2$ are the constants). \\

Since $0<-\Phi_4'(\nu)$ is increasing for $g_\nu\leq g_{\bar{\nu}}-A_1$ then

\bdis
-\Phi_4'(\nu)\geq \frac{1}{2\ln\frac{g_{\bar{\nu}}-A_1}{2\pi}}\ln\frac{2\pi mn}{g_{\bar{\nu}}-A_1}>0,\ mn>\frac{g_{\bar{\nu}}-A_1}{2\pi} ,
\edis

and, consequently,

\bdis
\bar{U}_{411}(r)=\mcal{O}\left(\frac{\ln T}{\ln\frac{mn}{Q_3}}\right),\ \frac{g_{\bar{\nu}}-A_1}{2\pi}=Q_3<mn\leq \frac T\pi\leq 2Q_3 ,
\edis

where $Q_3\in\mbb{N}$ (see Remark 4). Thus we obtain (comp. part. 10.2)

\be \label{10.14}
W_{43}=\mcal{O}\{(M+1)T^{1/2+\epsilon}\} ,
\ee

where $W_{43}$ is the contribution of $\bar{U}^1_{411}(r)$ into the sum $W_4$. \\

Since $0<\Phi_4'(\nu)$ is increasing for $g_\nu\geq g_{\bar{\nu}}+A_2$ then

\bdis
\Phi_4'(\nu)\geq \frac{1}{2\ln\frac{g_{\bar{\nu}}+A_2}{2\pi}}\ln\frac{g_{\bar{\nu}}+A_2}{2\pi mn}>0,\ mn<\frac{g_{\bar{\nu}}+A_2}{2\pi} ,
\edis

and consequently

\bdis
\bar{U}^2_{411}(r)=\mcal{O}\left(\frac{\ln T}{\ln\frac{Q_4}{mn}}\right),\ mn<Q_4=\frac{g_{\bar{\nu}}+A_2}{2\pi}
\edis

where $Q_4\in\mbb{N}$. Thus, we obtain (comp. the part 10.2)

\be\label{10.15}
W_{44}=\mcal{O}\{(M+1)T^{1/2+4\epsilon}\} ,
\ee

where $W_{44}$ is the contribution of $\bar{U}^2_{411}$ into the sum $W_4$. \\

The contribution of the term $\mcal{O}(\ln T)$ into $W_4$ is give by

\be \label{10.16}
\begin{split}
& W_{45}=\mcal{O}\left\{(M+1)\ln T\ssum_{T/2\pi\leq mn\leq T/\pi}\frac{d(m)d(n)}{\sqrt{mn}}\right\}= \\
& = \mcal{O}\{(M+1)\ln T\frac{T^{2\epsilon}}{\sqrt{T}}T^\epsilon T\}=\mcal{O}\{(M+1)T^{1/2+4\epsilon}\} .
\end{split}
\ee

\subsubsection{}

If $\tau_1<g_\nu\leq \frac{T}{\pi}$ then, similarly to the case $W_{41}$ (see the part 10.1) we obtain

\be \label{10.17}
W_{46}=\mcal{O}\{(M+1)T^{1/2+4\epsilon}\} .
\ee

Thus, from (10.1) by (10.7), (10.11), (10.14)-(10.17) the assertion (8.7) of the Lemma 6 follows.

\subsubsection{}

Finally, we complete the proof of the Lemma B. First of all, from (8.2) by (8.3), (8.4), (8.6), (8.7) we obtain

\be \label{10.18}
\bar{S}=\mcal{O}\{(M+1)T\ln^4T\} .
\ee

Hence, from (8.1) by (8.5), (10.18) the assertion (3.3) of the Lemma B follows.

\section{Some new classes of the formulae generated by the Theorem 1}

\subsection{}

For example, to the Euler's series

\bdis
\zeta(2)=\frac{\pi^2}{6}=\sum_{n=1}^\infty\frac{1}{n^2}=\sum_{n=1}^N\frac{1}{n^2}+\mcal{O}\left(\frac 1M\right)
\edis
where $N^2\leq M+1$ corresponds the class of formulae

\bdis
\sum_{T\leq t_\nu\leq 2T}Z^2(t_\nu)\sum_{n=1}^N Z^2\{ t_\nu\pm n\rho_1(\nu)\}\sim \frac{1}{8\pi^3}T\ln^5T,\ T\to\infty
\edis

for every random distribution of the signs $\pm$.

\subsection{}

In the connection with the discrete analog of the Hardy-Littlewood effect which we have proved in \cite{8},

\bdis
\frac{1}{N_1}\sum_{T\leq g_\nu\leq T+U}\left\{ \frac{1}{\bar{M}+1}\sum_{n=0}^{\bar{M}}Z(g_\nu+n\omega)\right\}^2\leq A\frac{\ln T}{\bar{M}},\
T\to\infty
\edis

where

\bdis
N_1=\sum_{T\leq g_\nu\leq T+U} 1,\ U=T^{5/12}\psi\ln^3T,\ \ln T<\bar{M}<\sqrt[3]{\psi}\ln T,\ \omega=\frac{\pi}{\ln\frac{T}{2\pi}} ,
\edis

we obtain from the formula (2.4) a more complicated - biquadratic - analog of the Hardy-Littlewood effect

\bdis
\begin{split}
&
\frac{1}{N_2}\sum_{T\leq g_\nu\leq 2T}\left\{ \frac{1}{M}\sum_{n=1}^{M}Z^2(t_\nu+n\rho_1(\nu))\right\}^2\sim
\frac 1\pi\frac{\ln^4T}{M},\ T\to\infty , \\
& N_2=\sum_{T\leq g_\nu\leq 2T} 1 .
\end{split}
\edis

\subsection{}

Let $\vp_n,\ n=1,2,\dots ,[t_i]$ be mutually independent random variables uniformly distributed within the segment $[-\pi,\pi]$. Next, let
$Z^2_\vp(t)$ be the random process generated by the phase-modulation with the random vector $(\vp_1,\dots ,\vp_{[t_i]})$ of the main term in the
Hardy-Littlewood's formula for $Z^2(t)$

\be \label{11.1}
Z^2_\vp(t)=2\sum_{n\leq[t_1]}\frac{d(n)}{\sqrt{n}}\cos\{ 2\vt(t)-t\ln n+\vp_n\},\ t_1=\frac{t}{2\pi} .
\ee

In this case, the formula (2.4) (and all its consequences) is valid for the class of all realizations of $Z^2_{\bar{\vp}}$ of the random process
(11.1). \\

\thanks{I would like to thank Michal Demetrian for helping me with the electronic version of this work.}

\end{document}